\documentclass[12pt]{article}
\usepackage{amsmath,amssymb,amsfonts}

\newtheorem{theorem}{Theorem}
\newtheorem{proposition}{Proposition}

\def\D{{\cal D}}
\def\A{{\cal A}}
\def\B{{\cal B}}
\def\R{{\mathbb R}}
\def\C{{\mathbb C}}

\begin{document}

\title{Blowing up solutions of the modified Novikov--Veselov equation and minimal surfaces}
\author{Iskander A. TAIMANOV
\thanks{Sobolev Institute of Mathematics, Academician Koptyug avenue 4, 630090, Novosibirsk, Russia, and Department of Mathematics and Mechanics, Novosibirsk State University, Pirogov street 2, 630090 Novosibirsk, Russia; e-mail: taimanov@math.nsc.ru. \newline
The work was supported by RSF (grant 14-11-00441). }}
\date{}
\maketitle

\section{Introduction}

In the present article we construct a solution to the modified Novikov--Veselov equation
(the two-dimensional generalization of the modified Korteweg--de Vries equation)
which has a singularity exactly at one point (Theorem 2).

The solution is given by an explicit formula
$$
\widetilde{U}(x,y,t) =
-\frac{3((x^2+y^2+3)(x^2-y^2)-6x(C-t))}{Q(x,y,t)},
$$
\begin{equation}
\label{u1}
Q(x,y,t) = (x^2+y^2)^3 + 3(x^4+y^4)+18 x^2 y^2 +9(x^2+y^2) +
\end{equation}
$$
+ 9(C-t)^2 +
(6x^3-18xy^2-18x)(C-t),
$$
from which it is clear that

\begin{itemize}
\item
{\sl it is infinitely differentiable (and even really-analytical) everywhere outside a single point
$x=y=0, t=C = \mathrm{const}$ at which it is not defined and has different finite limit values
along the rays $x/y = \mathrm{const}, t=C$, going into this point;}

\item
{\sl its restrictions
onto all planes $t=\mathrm{const}$ decay as $O(1/r^2)$, and, in particular, have finite $L_2$-norms;}

\item
{\sl the first integral (conservation law) $\int_{\R^2} \widetilde{U}^2 dx\, dy$ has the same value
equal to $3\pi$ for all times $t \neq C$ and jumps to $2\pi$ for $t=C$.}

\end{itemize}

The method of constructing such solutions is given by Theorem 1 and we consider in detail only only simplest example.
It is based on the geometrical interpretation \cite{T1} of the Moutard transformation for two-dimensional Dirac operators \cite{C}.

\section{Preliminary facts}

\subsection{The modified Novikov--Veselov equation}

The modified Novikov--Veselov (mNV) equation has the form
\begin{equation}
\label{mnv}
U_t = \big(U_{zzz} + 3U_z V + \frac{3}{2}UV_z \big) + \big(U_{\bar{z}\bar{z}\bar{z}} + 3U_{\bar{z}}\bar{V} + \frac{3}{2} U\bar{V}_{\bar{z}}\big),
\end{equation}
where
$$
V_{\bar{z}} = (U^2)_z,
$$
$z = x+iy \in \C$, $U$ is a real-valued function.

For making the equation correctly--posed we have to uniquely resolve the constraint which defines $V$.
For instance, for fast decaying solutions $U$ we may do that by assuming that $V$ is also fast decaying.

This equation takes the form of Manakov's $L,A,B$--triple:
$$
\D_t + [\D,\A] - \B\D = 0,
$$
where $\D$ is a two-dimensional Dirac operator:
\begin{equation}
\label{dirac}
\D = \left(
\begin{array}{cc}
0 & \partial \\
-\bar{\partial} & 0
\end{array}
\right) + \left(
\begin{array}{cc}
U & 0 \\
0 & U
\end{array}
\right),
\end{equation}
$\partial = \frac{\partial}{\partial z}$
and
$\bar{\partial} =
\frac{\partial}{\partial \bar{z}}$,
$$
\A = \partial^3 + \bar{\partial}^3 +
$$
\begin{equation}
\label{a}
+ 3\left(\begin{array}{cc}  V & 0 \\ U_z & 0 \end{array}\right)\partial +
3\left(\begin{array}{cc} 0 & -U_{\bar{z} }\\ 0 & \bar{V} \end{array}\right) \bar{\partial}+
\frac{3}{2}\left(\begin{array}{cc} V_z & 2U\bar{V} \\ -2UV & \bar{V}_{\bar{z}} \end{array}\right),
\end{equation}
$$
\B =
3\left(\begin{array}{cc} -V & 0 \\ -2U_z & V \end{array}\right)\partial +
3\left(\begin{array}{cc} \bar{V} & 2U_{\bar{z}} \\ 0 & -\bar{V} \end{array}\right)\bar{\partial} +
\frac{3}{2}\left( \begin{array}{cc} \bar{V}_{\bar{z}} - V_z & 2 U_{\bar{z}\bar{z}} \\ -2U_{zz} & V_z - \bar{V}_{\bar{z}}\end{array} \right)
$$

If $U$ depends only on $x$ and therewith $V = U^2$, then this equation reduces to the modified Korteweg--de Vries equation
$$
U_ t = \frac{1}{4} U_{xxx}  + 6 U_x U^2.
$$

The mNV equation was introduced in \cite{Bogdanov} and its name is due to the Novikov--Veselov equation introduced in \cite{NV,NV2})
which is a similar $2$-dimensional generalization of the Korteweg--de Vries equation.

\subsection{The Weierstrass representation of minimal surfaces}

A surface in $\R^3$ is called minimal if its mean curvature vanishes everywhere:
$$
H=0.
$$

The Weierstrass representation corresponds to every pair of holomorphic functions
$$
\psi_1, \bar{\psi}_2: {\cal U} \to \C
$$
a minimal surface
$$
F: {\cal U} \to \R^3,  \ \ \ \ F = (u^1,u^2,u^3),
$$
given by the formulas
$$
u^1(P) = \frac{i}{2} \int_{P_0}^P
\left((\psi_1^2 + \bar{\psi}_2^2)dz - (\bar{\psi}_1^2 + \psi_2^2)d\bar{z}\right) + u^1(P_0),
$$
\begin{equation}
\label{w}
u^2(P) = \frac{1}{2} \int_{P_0}^P
\left((\bar{\psi}_2^2 - \psi_1^2)dz + (\psi_2^2 - \bar{\psi}_1^2)d\bar{z}\right) + u^2(P_0),
\end{equation}
$$
u^3(P) = \int_{P_0}^P
\left(\psi_1 \bar{\psi}_2 dz + \bar{\psi}_1 \psi_2 \bar{z}\right) + u^3(P_0),
$$
where  $(u^1,u^2,u^3)$ are the Euclidean coordinates in $\R^3$,
$P_0 \in {\cal U}$ and the integral is taken over a path in ${\cal U}$ joining $P_0$ and $P$.
If ${\cal U}$ is simply-connected, then the integral does not depend on a choice of a path.
This is an immersion outside branch points where
the induced metric
\begin{equation}
\label{induced}
ds^2 = e^{2\alpha} dz d\bar{z} = (|\psi_1|^2 + |\psi_2|^2)^2 dz d\bar{z}
\end{equation}
vanishes.
The unit normal vector is equal to
$$
{\bf n} = \frac{1}{|\psi_1|^2 + |\psi_2|^2} ( i(\psi_1\psi_2 -\bar{\psi}_1\bar{\psi}_2),
- (\psi_1\psi_2 +\bar{\psi}_1\bar{\psi}_2),
(|\psi_2|^2 - |\psi_1|^2)).
$$
The formulas (\ref{w}) define a surface up to translations, i.e. up to $F(P_0)$.

It is well-known that every minimal surface in $\R^3$ admits such a representation.

\subsection{The Enneper surface}

The Enneper surface is an immersed (not embedded) minimal surface defined via the formulas
(\ref{w}) by
$$
\psi_1 = z, \ \ \ \psi_2 = 1.
$$
Substituting that into (\ref{w}), we obtain
$$
u^1(x,y) = y \left(\frac{y^2}{3}-x^2 -1\right) + u^1_0,
$$
\begin{equation}
\label{enneper}
u^2(x,y) = x \left(1+ y^2 - \frac{x^2}{3}\right) + u^2_0,
\end{equation}
$$
u^3(x,y) = x^2-y^2 + u^3_0,
$$
where ${\bf u}_0 = (u^1_0,u^2_0,u^3_0)$ is the image of the origin $x=y=0$
under an immersion.

\section{The Moutard transformation}

Let
$$
\psi = \left( \begin{array}{c} \psi_1 \\ \psi_2 \end{array}\right)
$$
be a solution of the Dirac equation
$$
\D \psi = 0
$$
where $\D$ is the Dirac operator (\ref{dirac}).
It is clear that
$$
\psi^\ast = \left(\begin{array}{c} -\bar{\psi}_2 \\ \bar{\psi}_1 \end{array}\right)
$$
satisfies the same equation. Let us
form a matrix-valued function $\Psi$ from $\psi$ and $\psi^\ast$ as follows
$$
\Psi =
\left(\begin{array}{cc}  \psi_1 & -\bar{\psi}_2 \\  \psi_2 & \bar{\psi}_1 \end{array}\right).
$$
It meets the matrix Dirac equation
\begin{equation}
\label{matrixeq}
\D \Psi = 0.
\end{equation}

We denote by $H$ a space formed by all matrices of the form
$$
\left(\begin{array}{cc} \alpha & \beta \\ -\bar{\beta} & \bar{\alpha} \end{array}\right), \ \ \ \alpha,
\beta \in \C,
$$
and put
$$
\Gamma =
\left(\begin{array}{cc} 0 & 1 \\ -1 & 0 \end{array}\right).
$$
It is evident that $H$ is closed under products and $\Gamma, \Psi \in H$.

For $U=0$ we have the operator
$$
\D_0  =
\left( \begin{array}{cc} 0 & \partial \\ -\bar{\partial} & 0
\end{array} \right)
$$
and vector functions $\psi$ which define minimal surfaces via (\ref{w}) are exactly solutions of
$\D_0\psi = 0$.

Given scalar functions $U$ and $V$, let us correspond to $H$-valued functions $\Phi$ and $\Psi$
 a  matrix-valued $1$-form
$$
\widetilde{\omega}(\Phi,\Psi) = \Phi^\top \Psi dy - i \Phi^\top \sigma_3 \Psi dx +
\left[i(\Phi^\top_{yy}\sigma_3\Psi + \Phi^\top\sigma_3\Psi_{yy}-\Phi^\top_y\sigma_3\Psi_y)
\right. +
$$
\begin{equation}
\label{omega1}
\left.
2iU(\Phi^\top_y\sigma_2\Psi - \Phi^\top\sigma_2\Psi_y)
+
\Phi^\top
\left(\begin{array}{cc}
iU^2-3iV & -iU_x \\
-iU_x & -iU^2 +3i\bar{V}
\end{array}\right)
\Psi
\right]dt =
\end{equation}
$$
-\frac{i}{2}\left(\Phi^\top \sigma_3 \Psi + \Phi^\top \Psi\right) dz - \frac{i}{2}\left(\Phi^\top \sigma_3 \Psi - \Phi^\top \Psi\right) d\bar{z} +
$$
$$
\left[-i((\Phi^\top_{zz} + \Phi^\top_{\bar{z}\bar{z}} - 2\Phi^\top_{z\bar{z}})\sigma_3\Psi +
\Phi^\top\sigma_3(\Psi_{zz} + \Psi_{\bar{z}\bar{z}} - 2 \Psi_{z\bar{z}}) -
\right.
$$
$$
(\Phi^\top_z - \Phi^\top_{\bar{z}})\sigma_3 (\Psi_z - \Psi_{\bar{z}})) -
2U((\Phi^\top_z - \Phi^\top_{\bar{z}}) \sigma_2 \Psi - \Phi^\top\sigma_2(\Psi_z - \Psi_{\bar{z}}))
+
$$
$$
\left.
\Phi^\top
\left(\begin{array}{cc}
iU^2-3iV & -i(U_z + U_{\bar{z}}) \\
-i(U_z + U_{\bar{z}}) & -iU^2 +3i\bar{V}
\end{array}\right)
\Psi
\right]dt,
$$
and matrix-valued functions
$$
\widetilde{S}(\Phi,\Psi)(z,\bar{z},t) = \Gamma \int_0^z \widetilde{\omega}(\Phi,\Psi),
$$
$$
K(\Psi) =  \Psi \widetilde{S}^{-1}(\Psi,\Psi)\Gamma \Psi^\top\Gamma^{-1}, \ \ \
$$
$$
M(\Psi) = \Gamma \Psi_y \Psi^{-1} \Gamma^{-1} =i\Gamma (\Psi_z - \Psi_{\bar{z}})\Psi^{-1}\Gamma^{-1}.
$$

The following Moutard transformation of solutions to the mNV equation was introduced
in \cite{C}.

\begin{proposition}
[\cite{C}]
\label{prop1}
Let $U(z,\bar{z},t)$ and $V(z,\bar{z},t)$ satisfy the mNV equation (\ref{mnv}), $\D$ is the family of Dirac operators with
potentials $U(z,\bar{z},t)$, and $\Psi_0(z,\bar{z},t)$ satisfy the system
$$
\D \Psi_0 = 0, \ \ \ \ \frac{\partial \Psi_0}{\partial t } = \A \Psi_0,
$$
where $\A$ has the form (\ref{a}).
Then

\begin{enumerate}
\item
the matrices $K(\Psi_0)$ and $M(\Psi_0)$ take the form
  $$
K =
\left(\begin{array}{cc} iW & a \\ -\bar{a} & -iW
\end{array}\right), \ \ \
M =
\left(\begin{array}{cc} b & c \\ -\bar{c} & \bar{b}
\end{array}\right),
$$
with $W$ real valued;

\item
for every solution $\Psi$ of the equations (\ref{matrixeq}) and
$$
\frac{\partial \Psi}{\partial t}= \A \Psi
$$
the function $\widetilde{\Psi}$ of the form
$$
\widetilde{\Psi} =  \Psi - \Psi_0 \widetilde{S}^{-1}(\Psi_0,\Psi_0) \widetilde{S}(\Psi_0,\Psi)
$$
satisfies the equations
$$
\widetilde{\D}\widetilde{\Psi} = 0
$$
for the Dirac operator $\widetilde{\D}$ with potential
\begin{equation}
\label{potential}
\widetilde{U} =  U + W
\end{equation}
and
$$
\frac{\partial \widetilde{\Psi}}{\partial t} = \widetilde{\A} \widetilde{\Psi}
$$
where $\widetilde{\A}$ takes the form (\ref{a}) with $U$ replaced by $\widetilde{U}$ and $V$ replaced by
$\widetilde{V}$:
$$
\widetilde{V} =  V + 2UW  + a^2 + 2(a\bar{b} - i\bar{c}W);
$$

\item
the function $\widetilde{U}$ is real-valued and $\widetilde{U}$ and $\widetilde{V}$ satisfy the mNV equation
\begin{equation}
\label{mnvt}
\widetilde{U}_t = \big(\widetilde{U}_{zzz} + 3\widetilde{U}_z \widetilde{V} +
\frac{3}{2}\widetilde{U}\widetilde{V}_z \big) + \big(\widetilde{U}_{\bar{z}\bar{z}\bar{z}} +
3\widetilde{U}_{\bar{z}}\bar{\widetilde{V}} + \frac{3}{2} \widetilde{U}\bar{\widetilde{V}}_{\bar{z}}\big),
\end{equation}
$$
\widetilde{V}_{\bar{z}} = (\widetilde{U}^2)_z
$$
\end{enumerate}
\end{proposition}

\section{Minimal surfaces and blowing up solutions of the mNV equation}
\label{section4}

Let us apply Proposition \ref{prop1} to the operator with $U=0$.
Although  this is a stationary solution of the mNV equation, the Moutard transformation leads to a  non-trivial non-stationary
solution of the mNV equation. A similar effect was found and used for the Novikov--Veselov equation
\cite{TT1,TT2}.

By straightforward computations we derive

\begin{theorem}
\label{th1}
Let $\psi_1(z,\bar{z},t)$ and $\psi_2(z,\bar{z},t)$ be a functions which satisfy the equations
$$
\bar{\partial}\psi_1 = \bar{\partial}\bar{\psi}_2 = 0,
$$
$$
\frac{\partial \psi_1}{\partial t}  = \frac{\partial^3 \psi_1}{\partial z^3}, \ \ \ \
\frac{\partial \psi_2}{\partial t}  = \frac{\partial^3 \psi_2}{\partial \bar{z}^3}.
$$
Then
\begin{equation}
\label{deforms}
\widetilde{S}(\Psi_0,\Psi_0)(z,\bar{z},t) =
\left(\begin{array}{cc}
iu^3 & -u^1 - iu^2 \\ u^1-iu^2 & -iu^3
\end{array}\right) +
i \int_0^t
\left(\begin{array}{cc}
w & \bar{v} \\ v & -w
\end{array}\right)d\tau,
\end{equation}
where
$$
\Psi_0 =
\left(\begin{array}{cc}  \psi_1 & -\bar{\psi}_2 \\  \psi_2 & \bar{\psi}_1 \end{array}\right),
$$
the minimal surfaces $F(z,\bar{z},t) = (u^1,u^2,u^3)$ are defined by $\psi_1$ and $\psi_2$ by
(\ref{w}) with ${\bf u}_0$ independent on $t$ and
$$
v = (\psi_{1,z}^2 - \psi_{2,\bar{z}}^2) - 2(\psi_1 \psi_{1,zz} - \psi_2\psi_{2,\bar{z}\bar{z}}),
$$
$$
w = \psi_{1,z}\bar{\psi}_{2,z} + \bar{\psi}_{1,\bar{z}}\psi_{2,\bar{z}}  - \psi_{1,zz}\bar{\psi}_2 - \psi_1 \bar{\psi}_{2,zz} -
\bar{\psi}_{1,\bar{z}\bar{z}} \psi_2 - \bar{\psi}_1 \psi_{2,\bar{z}\bar{z}}.
$$
 \end{theorem}

Theorem \ref{th1}
shows that $\widetilde{S}$ is a deformed minimal surface which depends on $t$
and is given by the second term in (\ref{deforms}).
We come to the following conclusion

\begin{itemize}
\item
{\sl to obtain a blowing up solution of the mNV equation we have to find a pair of $\psi_1$ and $\psi_2$ which satisfy the conditions of Theorem \ref{th1}
and such that the matrix $\widetilde{S}$ degenerates at some moment of time.}
\end{itemize}

The simplest candidate is given by the Enneper surface. In this case
$$
\psi_1 = z, \ \ \ \psi_2 = 1, \ \ \ v = 1, \ \ \ w=0.
$$
We put the image of the origin to be
$$
u^1_0 = u^3_0 = 0, \ \ \ u^2_0 = C > 0,
$$
and, by (\ref{enneper}), compute
$$
\widetilde{S}(x,y,t) =
\left(\begin{array}{cc}
iu^3 & -u^1 - iu^2 +it \\ u^1-iu^2 + it & -iu^3
\end{array}\right)
=
\left(
\begin{array}{cc} \gamma & \delta \\ -\bar{\delta} & \bar{\gamma}
\end{array}
\right)
$$
with
\begin{equation}
\label{constants}
\gamma = i(x^2-y^2), \ \ \
\delta = -y\left(\frac{y^2}{3}-x^2-1\right) -i\left[x\left(1+y^2-\frac{x^2}{3}\right)+C-t\right].
\end{equation}
We also easily derive that
$$
K =
\left(
\begin{array}{cc} z & -1 \\ 1 & \bar{z}\end{array}
\right)
\widetilde{S}^{-1}
\left(
\begin{array}{cc} \bar{z} & 1 \\ -1 & z\end{array}
\right),
$$
$$
M =
\frac{i}{1+|z|^2}\left(
\begin{array}{cc} -z & -1 \\ -1 & \bar{z}\end{array}
\right),
$$
and therefore
$$
W = -i\frac{|z|^2\bar{\gamma} + \gamma + \delta z - \bar{\delta}\bar{z}}{|\gamma|^2+|\delta|^2},
\ \ \ \
a = \frac{z(\bar{\gamma}-\gamma)-\delta z^2 -\bar{\delta}}{|\gamma|^2+|\delta|^2},
$$
$$
b = -\frac{iz}{1+|z|^2}, \ \ \
c = -\frac{i}{1+|z|^2}.
$$

Since $U=V=0$, we finally obtain
\begin{equation}
\label{u}
\widetilde{U} = -i\frac{|z|^2\bar{\gamma} + \gamma + \delta z - \bar{\delta}\bar{z}}{|\gamma|^2+|\delta|^2},
\end{equation}
\begin{equation}
\label{v}
\widetilde{V} =
\frac{(z(\bar{\gamma}-\gamma)-\delta z^2 -\bar{\delta})^2}{(|\gamma|^2+|\delta|^2)^2} +
\frac{2\widetilde{U}}{1+|z|^2}  - 2 \frac{iz(z(\bar{\gamma}-\gamma)-\delta z^2 -
\bar{\delta})}{(|\gamma|^2+|\delta|^2)(1+|z|^2)}
\end{equation}
where $\gamma$ and $\delta$ are given by (\ref{constants}).

Put $r = \sqrt{x^2+ y^2} = |z|$. It is clear that
$$
b = O\left(\frac{1}{r}\right), \ \ \ c = O\left(\frac{1}{r^2}\right) \ \ \ \ \mbox{as $r \to \infty$.}
$$
By (\ref{constants}), we have
$$
\gamma = O(r^2), \ \ \ \delta = O(r^3), \ \ \ a = O\left(\frac{1}{r}\right),
$$
and finally we derive that
\begin{equation}
\label{functions}
\widetilde{U}  = O\left(\frac{1}{r^2}\right), \ \ \
\widetilde{V}  = O\left(\frac{1}{r^2}\right) \ \ \ \ \mbox{as $r \to \infty$.}
\end{equation}

These functions $\widetilde{U}$ and $\widetilde{V}$ may have singularities only
at points where $|\gamma|^2+|\delta|^2 = 0$, i.e. exactly at the points where
the moving Enneper surface $(u^1, u^2-t, u^3)$ hits the origin. This motion preserves
$u^1$ and $u^3$ and, since we assume that $u^1_0=u^3_0=0$,
it is clear from (\ref{enneper}) that $u^1=u^3=0$ if and only if $x=y=0$.
However at $x=y=0$ we have $u^2 = C = \mathrm{const}$ and hence
$$
|\gamma|^2 + |\delta|^2 = 0 \ \ \ \mbox{if and only if} \ \ \ t=C.
$$

\begin{theorem}
\label{th2}
The functions $\widetilde{U}$ (\ref{u}) and $\widetilde{V}$ (\ref{v})
with $\gamma$ and $\delta$ given by (\ref{constants})

\begin{enumerate}
\item
satisfy the modified Novikov--Veselov equation (\ref{mnvt});

\item
decay at least quadratically in $r$: $\widetilde{U} = O(r^{-2}), \widetilde{V}=O(r^{-2})$;

\item
are really analytical $t \neq C$;

\item
have singularities exactly at $x=y=0, t=C$. At this
point $\widetilde{U}$ is not defined and
\begin{equation}
\label{sing}
\lim_{r \to 0, \varphi= \mathrm{const}} \widetilde{U}(z,\bar{z},C) = - \cos\,2\varphi \ \ \ \ \mbox{for $z = r e^{i\varphi}$};
\end{equation}

\item
\begin{equation}
\label{willmore}
\int_{\R^2} |\widetilde{U}|^2 dx \, dy =
\begin{cases} 3 \pi & \mbox{for $t \neq C$}, \\
2\pi & \mbox{for $t=C$}.
\end{cases}
\end{equation}
\end{enumerate}
\end{theorem}

The statements 1--3 of Theorem are established above.

The statement 4 follows from the formula (\ref{u1})
which is straightforwardly derived from (\ref{u}). For $C=t$ this formula reduces to
\begin{equation}
\label{u2}
\widetilde{U} = -\frac{3r^2(r^2+3)\cos\, 2\varphi}{r^2(r^4 + 3r^2 (1+\sin^2\,2\varphi)+9)}.
\end{equation}
Here $r$ and $\varphi$ are the polar coordinates:
$x+iy = ze^{i\varphi} = r(\cos\, \varphi + i \sin\,\varphi)$.

The statement 5 will be proved in the next section.

\section{Geometry of the blowing up solution of the mNV equation}

The exposition of the geometrical properties of the solution is based
on the explicit formulas for the action of the M\"obius inversion on the Weierstrass representation data
\cite{T1} and on the relation of $\int U^2 dxdy$ to the Willmore functional \cite{taimanov-mnv}.

Any surface in $\R^3$ is defined by the formulas (\ref{w})
(the Weierstrass representation) where a spinor $\psi$ satisfies the Dirac equation
$$
\D \psi = \left[\left(\begin{array}{cc} 0 & \partial \\ -\bar{\partial} & 0 \end{array}\right) +
\left(\begin{array}{cc} U & 0 \\ 0& U \end{array}\right)\right]\psi
$$
defines  a surface in $\R^3$ via formulas (\ref{w}).
Therewith $z$ is a conformal parameter on surface such that the induced metric takes the form
(\ref{induced}) and the real-valued potential $U$
is equal to
$$
U = \frac{e^\alpha H}{2} = \frac{(|\psi_1|^2+|\psi|^2)H}{2}
$$
with $H$ the mean curvature of the surface \cite{K,taimanov-mnv,T-RS}.
For $H=0$, i.e. for minimal surfaces, this representation is use the Weierstrass representation of minimal surfaces exposed above.

In \cite{T1} we show that

\begin{enumerate}
\item
the reduced matrix function
$$
S(\Psi_0,\Psi_0)(x,y,t) = \Gamma \int -\frac{i}{2}\left(\Psi_0^\top (\sigma_3+1)  \Psi_0 dz + \Psi_0^\top (\sigma_3 -1) \Psi_0 d\bar{z}\right),
$$
where
$$
\Psi_0 =
\left(\begin{array}{cc}
\psi_1 & -\bar{\psi}_2 \\ \psi_2 & \bar{\psi}_1
\end{array}
\right),
$$
is equal to
\begin{equation}
\label{s}
\left(\begin{array}{cc} iu^3 & -u^1 - iu^2 \\ u^1-iu^2 & -iu^3 \end{array}\right) \in su(2) \approx \R^3
\end{equation}
where $\Sigma_t = (u^1,u^2,u^3)$ is a surface defined up to translations by $\psi = \left(\begin{array}{c} \psi_1 \\
\psi_2 \end{array}\right)$  via (\ref{w}) at every moment $t$;

\item
the M\"obius inversion of $\R^3 \cup \{\infty\} = S^3$ in terms of (\ref{s}) takes a simple form
$$
S \to S^{-1}
$$
and if $\Psi_0$ defines a surface $\Sigma$, then the inverted surface is defined by
$$
\widetilde{\Psi} = \Psi_0 \cdot S^{-1}
$$
via the same formulas (\ref{w});

\item
the potential $U$ is transformed by the the inversion by the formula (\ref{potential}), i.e. by the Moutard transformation given in Proposition \ref{prop1}, with $\widetilde{S}$ replaced by $S$ in the definition of
$K(\Psi)$.
\end{enumerate}

The blowing up solution of the mNV equation exposed in Section \ref{section4}
has a very simple geometrical meaning:

\begin{itemize}

\item
the stationary function $\psi = \left(\begin{array}{c} z \\ 1 \end{array}\right)$
defines an immersion $S_0$ of the Enneper surface $\Sigma_0$, the matrix function $\widetilde{S}$
is equal to
$$
\widetilde{S}(x,y,t) = S_0(x,y) +
\left(\begin{array}{cc} 0 & it \\ it & 0 \end{array}\right)
$$
and defines a rigid translation $\Sigma_t$ of the initial Enneper surface along the $u^2$ axis: $u^2 \longrightarrow u^2-t$.

Since all surfaces $\Sigma_t$ are minimal, they have the same potential $U=0$, however the potentials of their
inversions $\Sigma_t^{-1}$ are different and are equal to  $\widetilde{U}(x,y,t)$ (\ref{u}).

The Enneper surfaces $\Sigma_t$ hit the origin only at one point $x=y=0$ and only at one moment of time $t=C$
and therewith the inversion maps this point into infinity, $\Sigma_C^{-1}$ becomes noncompact, and the potential
$\widetilde{U}$ achieves a singularity at $x=y=0$, $t=C$.

The quantity $4 \int \widetilde{U}^2 dx dy$ is the conservation law of the mNV equation and
is equal to the value of Willmore functional (the integral of the squared mean curvature) at the surface, i.e. in our case at $\Sigma^{-1}_t$
\cite{taimanov-mnv}.

\end{itemize}

Now the statement 5 of Theorem 2 follows, for instance, from computations of the values of Willmore functional for
inverted Enneper surfaces \cite{LN}.

\section{Final remarks}

1) The deformation $\Sigma^{-1}_t$ is an example of the mNV evolution of surfaces introduced
in \cite{K} for surfaces ``induced'' by the formulas (\ref{w}).

2) The constructed solution is special in many respects:

a) $\widetilde{S}(x,y,t)$ splits into $S_0(x,y) + P(t)$, i.e. describes a rigid motion of a minimal surface;

b) the inverted Enneper surfaces have many interesting geometrical features and, in particular, they are
branched Willmore spheres \cite{LN}.

Do rigid motions of other minimal surfaces in the same manner correspond to integrable soliton equations?

3) Other computable interesting examples can be found by using the higher order Enneper surfaces defined by the spinors $\psi = \left(\begin{array}{c} z^k \\ 1 \end{array}\right)$ and soliton spheres
(which are not minimal surfaces)
\cite{T-WN,BP}.

4) The results of this paper were briefly announced in \cite{TDAN}.

\end{document}